\documentclass[conference, onecolumn]{IEEEtran}
\usepackage{times}
\usepackage{amsmath}
\usepackage{amssymb}
\usepackage{latexsym}
\usepackage{graphicx, color}
\usepackage{epsfig}
\usepackage{url}
\usepackage{psfrag}
\usepackage{cite}
\usepackage{subfigure}

\newtheorem{theorem}{Theorem}
\newtheorem{lemma}{Lemma}

\newtheorem{observation}{Observation}

\newcommand{\remove}[1]{}

\begin{document}
\newcommand{\ls}[1]
   {\dimen0=\fontdimen6\the\font \lineskip=#1\dimen0
\advance\lineskip.5\fontdimen5\the\font \advance\lineskip-\dimen0
\lineskiplimit=.9\lineskip \baselineskip=\lineskip
\advance\baselineskip\dimen0 \normallineskip\lineskip
\normallineskiplimit\lineskiplimit \normalbaselineskip\baselineskip
\ignorespaces }

\title{\LARGE \bf
Dynamic Vehicle Routing for Data Gathering in Wireless Networks}

\author{\IEEEauthorblockN{G\"{u}ner D. \c{C}elik and Eytan Modiano}
\centering \IEEEauthorblockA{Massachusetts Institute of Technology
(MIT)\\ Email: \{gcelik, modiano\}@mit.edu}}

%


\maketitle \thispagestyle{empty} \pagestyle{plain}

\begin{abstract}
We consider a dynamic vehicle routing problem in wireless networks where messages arriving
randomly in time and space are collected by a mobile receiver
(vehicle or a collector). The collector is responsible for receiving these
messages via wireless communication by dynamically adjusting its
position in the network. Our goal is to utilize a combination of
\emph{wireless transmission} and \emph{controlled mobility} to
improve the delay performance in such networks. We show that the
necessary and sufficient condition for the stability of such a
system (in the bounded average number of messages sense) is given by
$\rho<1$ where $\rho$ is the average system load. We derive
fundamental lower bounds for the delay in the system and develop
policies that are stable for all loads $\rho<1$ and that have
asymptotically optimal delay scaling. Furthermore, we extend our
analysis to the case of multiple collectors in the network. We show
that the combination of mobility and wireless transmission results
in a delay scaling of $\Theta(\frac{1}{1-\rho})$ with the system
load $\rho$ in contrast to the $\Theta(\frac{1}{(1-\rho)^2})$ delay
scaling in the corresponding system where the collector visits each
message location.
\end{abstract}

%
\IEEEpeerreviewmaketitle

\section{Introduction}\ls{1.5}

\begin{figure}
\centering \psfrag{1\r}[l][][.5]{$\!\!\!\textrm{Collector}$}
\psfrag{2\r}[l][][.7]{$\Re$} \psfrag{3\r}[l][][.7]{$\!\!\!r^*$}
\psfrag{4\r}[l][][.5]{$\!\!\!\!\!\!\textrm{Message}$}
\psfrag{5\r}[l][][.7]{$\lambda$}
\includegraphics[width=0.3\textwidth]{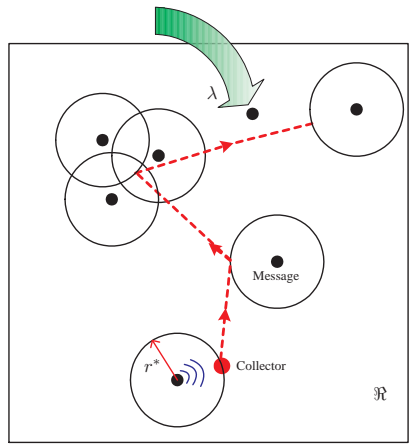}
\vspace{-1mm}\caption{The collector adjusts its position in order to
receive randomly arriving messages via wireless communication. The
circles with radius $r^*$ represent the communication range and the
dashed lines represent the collector's
path.}\label{Fig:intro_figure} \vspace{-5mm}
\end{figure}
There has been a significant amount of interest in performance
analysis of mobility assisted wireless networks in the last decade
(e.g., \cite{Tse}, \cite{Luo-Hub05}, \cite{Jerome08},\cite{emilio},
\cite{SH08}, \cite{Holly06},\cite{Holly07},\cite{ZhaMaYa08}).
Typically, throughput and delay performance of networks were
analyzed where nodes moving according to a random mobility model
were utilized for relaying data (e.g., \cite{GMPS04,Tse,NM05}). More
recently, networks deploying nodes with controlled mobility have been
considered focusing primarily on route design and ignoring the
communication aspect of the problem (e.g.,
\cite{Luo-Hub05}, \cite{Holly06}, \cite{Holly07},\cite{ZAZ2004}). %
%
In this paper we explore the use of \emph{controlled mobility} and
\emph{wireless transmission} in order to improve the delay
performance of wireless networks. We consider a dynamic vehicle
routing problem where a vehicle (collector) uses a combination of
physical movement and wireless reception to receive randomly
arriving data messages.

Our model consists of a collector that is responsible for gathering 
messages that arrive randomly in time at uniformly distributed
geographical locations. The messages are transmitted when the
collector is within their communication distance and depart the
system upon successful transmission. The collector adjusts its
position in order to successfully receive these messages in the
least amount of time as shown in Fig.~\ref{Fig:intro_figure}. This
setup is particularly applicable to networks deployed in a large
area so that a mobile element is necessary to provide connectivity
between spatially separated entities in the network. For instance,
this model is applicable to a sensor network where a mobile base
station collects data from a large number of sensors deployed at
random locations inside the network \cite{Burns-Levine-Inf05},
\cite{ZhaMaYa08}, \cite{ZAZ2004}. Another application is utilizing
Unmanned Aerial Vehicles (UAVs) as data harvesting devices or as
communication relays on a battlefield environment \cite{Jerome08}.
This model also applies to networks in which data rate is relatively
low so that data transmission time is comparable to the collector's
travel time, for instance in underwater sensor networks
\cite{Akyildiz_05}, \cite{Soma-kansal-06}.

Vehicle Routing Problems (VRPs) have been extensively
studied in the past (e.g., \cite{arkinGCSP94}, \cite{bert1},
\cite{bert2}, \cite{bert3}, \cite{FB2004}, \cite{MitTSPN07}, \cite{Holly06},
\cite{Holly07}, \cite{Xu}). The common example of a VRP is the
Euclidean Traveling Salesman Problem (TSP) in which a single server
is to visit each member of a fixed set of locations on the plane
such that the total travel cost is minimized. Several extensions of
TSP have been considered in the past such as stochastic demand
arrivals and the
use of multiple servers \cite{bert1}, \cite{bert2}, \cite{FB2004}. 
In particular, in the TSP
with neighborhoods (TSPN) problem 
the vehicle is to visit a neighborhood of each demand location
\cite{arkinGCSP94},
\cite{MitTSPN07}, which can model a mobile collector receiving
messages from a communication distance. A more detailed review of
the literature in this field can be found in \cite{bert2},
\cite{MitTSPN07} and \cite{Xu}. 

Of particular relevance to us among the VRPs is the Dynamic
Traveling Repairman Problem (DTRP) due to Bertsimas and van Ryzin
\cite{bert1}, \cite{bert2}, \cite{bert3}. 
DTRP is a stochastic and dynamic VRP in which a vehicle is to serve
demands that arrive randomly in time and space. 
Fundamental lower bounds on delay were established and several
vehicle routing policies were analyzed for DTRP for a single server
in \cite{bert1}, for multiple servers in \cite{bert2} and for
general demand and interarrival time distributions in \cite{bert3}.
Later, \cite{Holly06}, \cite{Holly07} generalized this model to
analyze Dynamic Pickup and Delivery Problem (DPDP) where fundamental
bounds on delay were established. We apply this model to wireless
networks where the demands are data messages to be transmitted to a
collector which is capable of wireless communication\footnote{In
\cite{bert1}, \cite{bert2}, or \cite{bert3} the collector needs to be at the message
location in order to be able to serve it, therefore, we will refer to the DTRP
model as the system without wireless transmission.}. In our system
the problem has considerably different characteristics since in this
case the collector does not have to visit message locations but
rather can receive the messages from a distance using wireless
communication. The objective in our system is to effectively utilize
this combination of wireless transmission and controlled mobility in
order to minimize the time average message waiting time. 

In a closely related problem where multiple mobile nodes with
controlled mobility and communication capability relay the messages
of static nodes, \cite{emilio} derived a lower bound on node travel
times. Message sources and destinations are modeled as static nodes
in \cite{emilio} and these nodes have saturated arrivals hence
queuing aspects were not considered. In an independent work, \cite{Ver-Altm-09} considered utilizing mobile wireless servers as data relays on periodic routes and
applied various delay relations from Polling models to this setup. A mobile server harvesting data
from spatial queues in a wireless network was considered in
\cite{Jerome08} where the stability region of the system was
characterized using a fluid model approximation. In \cite{Celik} we
analyzed a model similar to the current paper but for which the arriving messages were
transmitted to the collector using a random access scheme, creating
interference among neighboring transmissions. In this paper, the
message transmissions are scheduled, i.e., there is only one
transmission in the system at a given time, and the collector
decides on the message to be transmitted next. The two systems have considerably different
characteristics as will be explained in the following sections.

Another related body of literature lies in the area of utilizing
mobile elements that can control their mobility to collect
sensor data in Delay Tolerant Networks (DTN) 
(e.g., \cite{Burns-Levine-Inf05,Luo-Hub05,SH08,Soma-kansal-06,ZhaMaYa08,ZAZ2004}). 
Route selection (e.g., \cite{Luo-Hub05},
\cite{SH08},\cite{ZAZ2004}), scheduling or dynamic mobility control
(e.g., \cite{Burns-Levine-Inf05}, 
\cite{Soma-kansal-06},\cite{ZhaMaYa08}) algorithms were proposed to
maximize network lifetime, to provide
connectivity or to minimize delay. 
These works focus primarily on mobility and usually consider
particular policies for the mobile element. To the best of our
knowledge,
this is the first attempt to develop fundamental bounds on delay 
in a system where a collector is to gather data messages randomly
arriving in time and space using \emph{wireless communication} and
\emph{controlled mobility}.

The main contributions of this work are the following. We show that
$\rho<1$ is the necessary and sufficient condition for the stability
of the system where $\rho$ is the system load. We derive fundamental
lower bounds on delay and develop algorithms that
are asymptotically within a constant factor of the lower bounds. 
We show that the combination of mobility and wireless transmission
results in a delay scaling of 
$\Theta(1/(1-\rho))$ 
in contrast to the $\Theta(1/(1-\rho)^2)$ delay scaling in the
system where the collector visits each message location analyzed in
\cite{bert1}, \cite{bert2}. Finally, we extend our analysis
to the case of multiple collectors in the network. 

This paper is organized as follows. In 
Section \ref{Sec:Model} we describe the system model and in Section
\ref{Sec:Stab} we characterize the necessary and sufficient
conditions for the stability of the system. We derive fundamental
lower bounds on delay in Section \ref{Sec:LB_Delay} and in Section
\ref{Sec:Policies} we provide upper bounds on delay together with
numerical results.
We extend the analysis for a system with multiple collectors in
Section \ref{Sec:Mult_veh}.

%
%
\section{Model}\label{Sec:Model}

Consider a square region 
$\mathcal{R}$ of area $A$ and messages arriving into $\mathcal{R}$ according to a
Poisson process (in time) of intensity $\lambda$. Upon arrival the
messages are distributed independently and uniformly in $\mathcal{R}$ and
they are to be gathered by
a collector via wireless reception. 
An arriving message is transmitted to the collector when the
collector comes within the \emph{reception distance} of the message
location and grants access for the message's transmission. Therefore, there is no interference power from the neighboring nodes during message receptions. 
%
%
%
%
%
We assume that the transmit power $P_T$ is constant and that the
transmissions are subject to distance attenuation.
In such a system, the received power of a transmission from node $i$,
located at distance $r_i$ from the collector, is given by $P_{R,i} =
P_{T}Kr_i^{-\alpha}$ \cite{Tse}, \cite{GuptaKumar},
where 
$\alpha$ is the power loss exponent (typically between $2$ and $6$),
and $K$ is the attenuation constant normalized to $1$. 

Next we argue that the Signal to Noise Ratio (SNR) packet reception model \cite{Tse}, \cite{GuptaKumar} is equivalent to a disk model \cite{Celik}, \cite{GuptaKumar} under the above assumptions.
In the SNR model, a transmission is successfully decoded at the collector if its
SNR is above a threshold $\beta$, i.e., if
$\textrm{SNR}_i=P_{R,i}/{P_N} \ge \beta,$
where $P_N$ is the background noise power.
We let $r$ be the \emph{reception distance} of the collector.
If the location of the next message to be received is within 
$r$, the collector stops and attempts
to receive the message. Otherwise, the collector travels towards the
message location until it is within a distance $r$ away from the
message. 
A transmission at distance $r$ to the collector
is successful if
$r \le ( \textrm{SNR}_c/\beta )^{1/\alpha}$ where $\textrm{SNR}_c=P_T/P_N$ denotes the SNR of a transmission
from unit distance. Therefore, the optimal
reception distance is the maximum reliable communication distance
$r^* = (\textrm{SNR}_c/\beta)^{1/\alpha}$. Hence, essentially we have a disk model of radius $r^*$, where a transmission can be received only if it is within a disk of radius $r^*$ around the collector.
Under this model, transmissions are assumed to be at a constant rate
taking a fixed amount of time denoted by
$s$.


%

The collector travels from the current message reception point to
the next message reception point at a constant speed $v$. 
We assume that at a given time the collector knows the locations and
the arrival times of the messages that arrived before this time.
The knowledge of the service locations is a standard assumption in vehicle routing literature
\cite{arkinGCSP94}, \cite{bert1}, 
\cite{FB2004}, \cite{MitTSPN07}, \cite{Holly06}. 
Let $N(t)$ denote the total number of
messages in the system at time $t$. We say that the system is stable
under a policy if
\cite{Asmussen}, 
\cite{ModShahZuss06}, 
\begin{equation}
\displaystyle \limsup_{t\rightarrow \infty} \mathbb{E} [N(t)] <
\infty,
\end{equation}
namely, the long term expected number of messages in the system is
finite. Let $\rho=\lambda s$ denote the load arriving into the system per unit time. For
stable systems, $\rho$ denotes the fraction of time the
collector spends receiving messages.

We define $T_i$ as the time between the
arrival of message $i$ and its successful reception. $T_i$ has three
components: $W_{d,i}$, the waiting time due to collector's travel
distance from the time message $i$ arrives until it gets served,
$W_{s,i}$, the waiting time due to the reception times of messages
received from the time message $i$ arrives until it gets served, and
$s$, reception time of the message. 
The total waiting time of message $i$ is denoted by
$W_i=W_{d,i}+W_{s,i}$, hence $W_i=T_i-s$.  The fraction of time the
collector spends receiving messages is denoted by $\rho$, and for
stable systems $\rho=\lambda s$. We let $d_i$ be the collector
travel distance from the collector's reception location for the
message served prior to message $i$ to collector's reception
location for message $i$.
%
The time average per-message travel distance of the collector,
denoted by $\bar{d}$, is defined by an expectation in the steady
state given by $\bar{d}=\lim_{i
\rightarrow \infty} E[d_i]$. 
The time average delays $T$, $W$, $W_d$ and $W_s$ are defined
similarly to have $T=W_d+W_s+s$ where all the limits are assumed to
exist. $T^*$ is defined to be the optimal system time which
is given by the policy that minimizes $T$. 

%
%


\section{Stability}\label{Sec:Stab}
In this section we characterize a necessary and sufficient condition
for the stability of the system.

%
%

\subsection{Necessary Condition for Stability}
\begin{theorem}\label{thm:inf_vel_sin_veh}
\emph{A necessary condition for the stability of any policy is $\rho<1$.
Furthermore, the time average waiting time satisfies}
\begin{equation}\label{eq:delay_PK}
W \ge \frac{\lambda s^2}{2(1-\rho)}.
\end{equation}
\end{theorem}
\begin{proof}
We first show that the unfinished work and the delay experienced by
a message in the system stochastically dominates that in the
equivalent system with zero travel times for the collector.
\begin{lemma}\label{lem:app_pf}
\emph{The steady state time average delay in the system is at least
as big as the delay in the equivalent system in which travel times
are considered to be zero (i.e., $v=\infty$).}
\end{lemma}
\begin{proof}
Consider the summation of per-message reception and travel times,
$s$ and $d_i$, as the total service requirement of a message in each
system. Since $d_i$ is zero for all $i$ in the infinite velocity
system and since the reception times are constant equal to $s$ for
both systems, the total service requirement of each message in our
system is deterministically greater than that of the same message in
the infinite velocity system. Let $D_1,D_2,...,D_n$ and
$D^{'}_1,D^{'}_2,...,D^{'}_n$ be the departure instants of the
$1^{st}$, $2^{nd}$ and similarly the $n^{th}$ message in the
original and the infinite velocity system respectively. 
Similarly let $A_1,A_2,...,A_n$ be the arrival times of the
$1^{st}$, $2^{nd}$ and the $n^{th}$ message in both systems. We will
use induction to prove that $D_i \ge D_i^{'}$ for all $i$. Since the
service requirement of each message is smaller in the infinite
velocity system, we have $D_1 \ge D_1^{'}$. Assuming we have $D_n
\ge D_n^{'}$, we need to show that $D_{n+1} \ge D_{n+1}^{'}$.
\begin{equation}\label{eq:app_proof_temp1}
A_{n+1} \le D_{n+1}-s,
\end{equation}
hence the $n+1^{th}$ message is available before the time
$D_{n+1}-s$. We also have
\begin{equation*}
D_{n}^{'} \le D_{n} \le D_{n+1}-s.
\end{equation*}
The first inequality is due to 
the induction hypothesis and the second inequality is because we
need at least $s$ amount of time between the $n^{th}$ and $n+1^{th}$
transmissions. Hence the collector is available in
the infinite velocity system before the time $D_{n+1}-s$. Combining
this with (\ref{eq:app_proof_temp1}) proves the induction. Now let
$D(t)$ and $D^{'}(t)$ be the total number of departures by time $t$
in our system and the infinite velocity system respectively.
Similarly let $N(t)$ and $N^{'}(t)$ be the total number of messages
in the two systems at time $t$. Finally let $A(t)$ be the total
number of arrivals by time $t$ in both systems. We have
$A(t)=N(t)-D(t)=N^{'}(t)-D^{'}(t)$. From the above induction we have
$D(t)\le D^{'}(t)$ and therefore
\begin{equation*}
N(t) \ge N^{'}(t).
\end{equation*}
Since this is true at all times, we have that the time average
number of customers in the system is greater than that in the
infinite velocity system. Finally using Little's law proves the
lemma.
\end{proof}
Since the infinite velocity system behaves as an M/G/1
queue (an M/G/1 queue is a queue with Poisson arrivals, general
i.i.d. service times and 1 server and an M/D/1 queue has constant
service times), the average waiting time in this system is given by
the Pollaczek-Khinchin (P-K) formula for M/G/1 queues \cite[p.
189]{BertGall92}, i.e., $\lambda s^2/(2(1-\lambda s))$. Therefore we
have (\ref{eq:delay_PK}). Furthermore, a direct consequence of this
lemma is that a necessary condition for stability in the infinite
velocity system is also necessary for our system. It is well-known
that the necessary (and sufficient) condition for stability in the
M/G/1 systems is given by $\rho<1$ (see e.g., \cite{BertGall92} or
\cite{Geor_Neely_Tass06}). 
\end{proof}

\subsection{Sufficient Condition for Stability}
Here we prove that $\rho < 1$ is a sufficient condition for
stability of the system under a policy based on Euclidean TSP with
neighborhoods (TSPN). TSPN is a generalization of TSP in which the
server is to visit a neighborhood of each demand location via the
shortest path \cite{arkinGCSP94}, \cite{MitTSPN07}. In our case the
neighborhoods are disks of radius $r^*$ around each message
location. TSPN is an NP-Hard problem such as TSP. Recently,
\cite{MitTSPN07} proved that a Polynomial Time Approximation Scheme
(PTAS) exists for TSPN among fat regions in the plane. A region is
said to be \emph{fat} if it contains a disk whose size is within a
constant factor of the diameter of the region, e.g., a disk, and a
PTAS belongs to a family of $(1 + \epsilon)$-approximation
algorithms parameterized by $\epsilon > 0$.
%
%

\subsubsection{TSPN Policy}\label{Sec:TSPN_Pol}

Assume the system is initially empty (at time $t_0=0$). The receiver
waits at the center of $\mathcal{R}$ until the first message arrival, moves
to serve this message and returns to the center. Let time $t_1$ be
the time at which the receiver returns to the center. At $t_1$, if
the system is empty, the receiver repeats the above process and we
define $t_2$ similarly. If there are messages waiting for service at
time $t_1$, the receiver computes the TSPN tour (e.g., using the
PTAS in \cite{MitTSPN07}) through all the messages that are present
in the system at time $t_1$, receives these messages in that tour and returns to the center. 
We let $t_2 > t_1$ be the first time when the receiver returns to
the center after receiving all the messages that were present in the
system at $t_1$ and repeat the above process. 
We define the epochs $t_i$ as the time the receiver returns to the
center after serving all the messages that were present in the
system at time $t_{i-1}$ \footnote{A similar policy based on TSP was
discussed in \cite{PavFraz07} for a system without communication
capability similar to DTRP, where the time average delay of the
policy was characterized for the heavy load regime.}.

Let the total number of messages waiting for service at time $t_i$,
$N_i\triangleq N(t_i)$, be the system state at time $t_i$. 
Note that $N_i$ is an irreducible Markov chain on countable state
space $\textbf{N}$. We show the stability of the TSPN policy
through the ergodicity of this Markov chain. 

\begin{theorem}\label{thm:TSPN}
The system is stable under the TSPN policy for all loads $\rho <1$.
\end{theorem}
\begin{proof}
Given the system state $N_i$ at time $t_i$, we apply the algorithm
in \cite{MitTSPN07} to find a TSPN tour of length $L_i$ through the
$N_i$ neighborhoods that is at most $(1+\epsilon)$ away from the
optimal TSPN tour length $L^*_i$. Note that $L^*_i$ can be upper bounded
by a constant $L$ for all $N_i$. This is because the collector does not have to move for messages within its communication range and a finite number of such disks of radius $r^*$ can cover the network region for any $r^*>0$. The collector then can serve the messages in each disk from its center incurring a tour of constant length $L$ (an example of such a tour is shown in Fig. \ref{Fig:part_pol}). We will use the Lyapunov-Foster criterion to show that the Markov
chain described by the states $N_i$ is positive recurrent
\cite{Asmussen}. We use $V(N_i)=sN_i$, the total load served during
$i^{th}$ cycle, as the Lyapunov function (note that $V(0)=0$,
$S_k=\{x:V(x) \le K\}$ is a bounded set for all finite $K$ and
$V(.)$ is a non-decreasing function). Since the arrival process is
Poisson, the expected number of arrivals during a cycle can be
upper-bounded as follows:
\begin{equation}\label{eq:N_i_TSPN}
\mathbb{E} [N_{i+1}| N_i] \le \lambda (L/v
+ sN_i ).
\end{equation}
Hence we obtain the following drift expression for the load during a
cycle.
\begin{equation*}
\mathbb{E} [sN_{i+1}- sN_{i}| N_i] \le
\rho L/v - (1-\rho) sN_i.
\end{equation*}
Since $\rho<1$, there exist a $\delta>0$ such that $\rho+\delta <1$: 
\begin{eqnarray}
\mathbb{E} [sN_{i+1}- sN_{i}| N_i] &\le&
\rho L/v - \delta sN_i\nonumber\\
&\le& -\delta s  + \frac{\rho L}{v}.\mathbf{1}_{\{N_i\in
 S\}},
\end{eqnarray}
where $\mathbf{1}_{\{N\in S\}}$ is equal to $1$ if $N\in S$ and zero
otherwise and $S = \{N \in \textbf{N}: N \le K \}$ is a bounded set
with $K=\lceil \frac{\rho L }{v\delta s}+1\rceil$. Hence the drift is negative as long as $N_i$ is outside a bounded set.
Therefore, by the standard Lyapunov-Foster criterion
\cite{AltKonsLiu92}, \cite{Asmussen}, the Markov chain $(N_i)$ is
positive recurrent, it has a unique stationary distribution and we
can bound the steady state time average of $N_i$ as
\cite{ModShahZuss06}
\begin{equation}
\displaystyle \limsup_{t_i\rightarrow \infty} \mathbb{E} [N(t_i)]
\le \frac{\lambda L}{v(1-\rho)}.
\end{equation}
Furthermore, given some $t\in [t_i,t_{i+1}]$, 
we have
\begin{eqnarray}
\displaystyle \limsup_{t\rightarrow \infty} \mathbb{E} [N(t)] &&
\!\!\!\!\!\!\!\!\le \displaystyle \limsup_{t_i\rightarrow \infty}
\mathbb{E}
[N(t_i) + N(t_{i+1})] \nonumber \\
&&\!\!\!\!\!\!\!\!\le
2\frac{\lambda L }{v(1-\rho)} < \infty. \label{eq:TSPN_delay}
\end{eqnarray}
\end{proof}
The delay scaling of the TSPN policy with load $\rho$ is $\frac{1}{1-\rho}$ as shown in (\ref{eq:TSPN_delay}), the same delay scaling as in a G/G/1 queue. This is a fundamental improvement in delay due to the communication capability as the system without wireless transmission in \cite{bert1} has $\Theta(\frac{1}{(1-\rho)^2})$ delay scaling.

Note that $\rho<1$ is a sufficient stability condition also for the system without communication capability. This case corresponds to $r^*=0$, where we utilize a $(1+\epsilon)$ PTAS for the optimal TSP tour
through the message locations instead of the TSPN tour.
An upper bound on the TSP tour for any $N_i$ points arbitrarily
distributed in a square of area $A$ is given by
$\sqrt{2AN_i}+1.75\sqrt{A}$ \cite{Lawler85}. Similar arguments as above leads to the drift condition
%
%
\begin{equation*}
\mathbb{E} [sN_{i+1}- sN_{i}| N_i] \le
\rho(\kappa_1\sqrt{AN_i}+\kappa_2) - (1-\rho) sN_i,
\end{equation*}
for some constants $\kappa_1$ and $\kappa_2$, where
the drift is again negative as long as $N_i$ is outside a  bounded
set $S$. The difference in this case is that the travel time per cycle scales with the number of messages $N_i$ as $\sqrt{N_i}$  which can be shown to result in $O(\frac{1}{(1-\rho)^2})$ delay scaling with the load $\rho$.

\section{Lower Bound On Delay}\label{Sec:LB_Delay}
For wireless networks with a small area and/or very good channel
quality such that $r^* \ge \sqrt{A/2}$, the collector does not need
to move as every message will be in its reception range if it just
stays at the center of the network region. In that case the system
can be modeled as an M/D/1 queue with service time $s$
and the associated queuing delay is given by the 
P-K formula for M/G/1 queues, i.e., $W=\lambda s^2/(2(1-\rho))$. However, when $r^* < \sqrt{A/2}$, the collector has to move in
order to receive some of the messages. In this case the reception
time $s$ is still a constant, however, the travel time per message is now a random variable which is not independent over
messages (for example, observing small travel times for the previous messages implies a
dense network, and hence the future travel times per message are also
expected to be small). Next we provide a lower bound similar to a
lower bound in \cite{bert1} with the added complexity of
communication capability in our system. 
\begin{theorem}\label{thm:LB_delay}
\emph{The optimal steady state time average delay $T^*$ is lower
bounded by\footnote{Note that $(||U||-r^*)^+$ represents $\max ( 0,
||U||-r^*)$ and $U$ is a uniformly distributed random variable over
the network region $\mathcal{R}$.}
\begin{equation}\label{eq:T_LB_low_lam}
T^* \ge \frac{ \mathbb{E} [ (||U||-r^*)^+] }{v(1-\rho)} +
\frac{\lambda s^2}{2(1-\rho)} + s.
\end{equation}
}
\end{theorem}
\begin{proof}
As outlined in Section \ref{Sec:Model}, the delay of message $i$, $T_i$ has three components: $T_i=W_{d,i} +
W_{s,i} + s_i$.
Taking expectations and the limit as $i \rightarrow \infty$ yields
\begin{equation}\label{eq:low_lam_del2}
T=W_{d} + W_{s} + s.
\end{equation}
A lower bound on $W_d$ is found as follows: Note that $W_{d,i}.v$ is
the average \emph{distance} the collector moves during the waiting
time of message $i$. This distance is at least as large as the
average distance between the location of message $i$ and the
collector's location at the time of message $i$'s arrival less the
reception distance $r^*$. The location of an arrival is determined
according to the uniform distribution over the network region, while
the collector's location distribution is in general unknown as it
depends on the collector's policy. We can lower bound $W_d$ by
characterizing the expected distance between a uniform arrival and
the best a priori location in the network that minimizes the
expected distance to a uniform arrival. Namely we are after the
location $\nu$ that minimizes $\mathbb{E} [||U-\nu||]$ where $U$ is
a uniformly distributed random variable. The location $\nu$ that
solves this optimization is called the \emph{median} of the region
and in our case the median is the center of the square shaped
network region. Because the travel distance is nonnegative, we
obtain the following bound on $W_d$:
\begin{equation}\label{eq:Wd_temp}
W_d \ge \frac{\mathbb{E} [ (||U||-r^*)^+] }{v}.
\end{equation}
Let $N$ be the average number of messages received in a waiting time
and let $R$ be the average residual reception (service) time. Due to
the PASTA property of Poisson arrivals (Poisson Arrivals See Time
Averages) (see for example \cite[p. 171]{BertGall92}) a given
arrival in steady state observes the time average steady state
occupancy distribution. Therefore, the average residual time
observed by an arrival is also $R$ and is given by $\lambda
s^2/2$ \cite[p. 188]{BertGall92} and we have
\begin{equation}\label{eq:Ws_temp1}
W_s = sN + R.
\end{equation}
Since in a stable system in steady state the average number of
messages received in a waiting time is equal to the average number
of arrivals in a waiting time (a variation of Little's law) we have
$N= \lambda W= \lambda (W_d+W_s)$. Substituting this in
(\ref{eq:Ws_temp1}) we obtain
\begin{equation*}
W_s = s \lambda (W_d+W_s) + \frac{\lambda s^2}{2}.
\end{equation*}
This implies
\begin{equation}\label{eq:Ws_temp}
W_s = \frac{\rho}{1-\rho}W_d + \frac{\lambda
\overline{s^2}}{2(1-\rho)},
\end{equation}
Substituting (\ref{eq:Wd_temp}) and (\ref{eq:Ws_temp}) in
(\ref{eq:low_lam_del2}) 
yields (\ref{eq:T_LB_low_lam}).
\end{proof}
In addition to the average waiting time of a
classical M/G/1 queue given in (\ref{eq:delay_PK}), the queueing
delay also increases due to the
collector's travel. Note that the $\mathbb{E} [ (||U||-r^*)^+]$ term can be further lower bounded by
$\mathbb{E}[||U||]-r^*$, where $\mathbb{E}[||U||]=0.383\sqrt{A}$ \cite{bert1}.

\section{Collector Policies}\label{Sec:Policies}

We derive upper bounds on delay via analyzing policies for the
collector. The TSPN policy analyzed in Section \ref{Sec:TSPN_Pol} is
stable for all loads $\rho<1$ and has $O(\frac{1}{1-\rho})$ delay scaling. 
Since the lower bound in Section \ref{Sec:LB_Delay} also scales with the load as
$\frac{1}{1-\rho}$, \emph{the TSPN  policy has optimal delay scaling}. 

\subsection{First Come First Serve (FCFS) Policy}
A straightforward policy is the FCFS policy
where the messages are served in the order of their arrival times. A
version of the FCFS policy, call FCFS', where the receiver has to
return to the center of the network region (the median of the region
for general network regions) after each message reception is shown
to be optimal at light loads for the DTRP problem
\cite{bert1}. This is because the center of the network region is
the location that minimizes the
expected distance to a uniformly distributed arrival. 
Since in our system we can do at least as good as the DTRP, FCFS' is
optimal also for our system at light loads. Furthermore, the FCFS
policy is not stable for all loads $\rho <1$, namely, there exists a
value $\hat{\rho}$ such that the system is unstable under FCFS
policy for all $\rho>\hat{\rho}$. This is because in the FCFS system
the average travel component of the service time is fixed, which
makes the average arrival rate greater than the average service rate
as $\rho \rightarrow 1$. Therefore, it is better for a policy to
serve more messages in the same ``neighborhood'' in order to reduce
the amount of time spent on mobility.

\subsection{Partitioning Policy}\label{Sec:Part_pol}
Next we propose a policy based on partitioning the network region
into subregions and the collector performing a cyclic service of the
subregions. This policy is an adaptation of the Partitioning policy
of \cite{bert1} to the case of a system with wireless transmission.
We explicitly derive the delay expression for this policy and
show that it scales with the load as $O(\frac{1}{1-\rho})$ as in the TSPN policy.

\begin{figure}
\centering \psfrag{1\r}[l][][.7]{$\!r^*$}
\psfrag{2\r}[l][][.7]{$\!\!\!\!\sqrt{2}r^*$}
\includegraphics[width=0.25\textwidth]{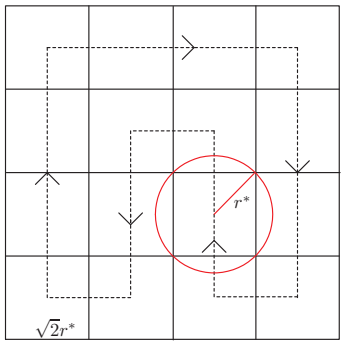}
\vspace{-1mm}\caption{The partitioning of the network region into
square subregions of side $\sqrt{2}r^*$. The circle with radius
$r^*$ represents the communication range and the dashed lines
represent the collector's path.} \label{Fig:part_pol} \vspace{-5mm}
\end{figure}

We divide the network region into
$(\sqrt{2}r^* \,\textrm{x}\, \sqrt{2}r^*)$ squares as shown in Fig.
\ref{Fig:part_pol}. 
This choice ensures us that every location in the square is within
the communication distance $r^*$ of the center of the square.
The number of subregions in such a Partitioning is given
by\footnote{If $\sqrt{n_s}$ is not even or if the subregions do not
fit the network region for a particular choice of $r^*,$ then one
can partition the region using the largest reception distance
$\underline{r}^* < r^*$ such that these conditions are satisfied.}
$n_s =A/(2(r^*)^2)$.
The partitioning in Fig. \ref{Fig:part_pol} represents the case of
$n_s=16$ subregions.
The collector services the subregions in a cyclic order as displayed
in Fig. \ref{Fig:part_pol} by receiving the messages in each
subregion from its center using an FCFS order. 
The messages within each subregion are served exhaustively, i.e.,
all the messages in a subregion are received before moving to the
next subregion. The collector then receives the messages in the next
subregion exhaustively using FCFS order and repeats this process. 
The distance traveled by the collector between each subregion is a
constant equal to $\sqrt{2}r^*$. It is easy to verify that the
Partitioning policy behaves as a multiuser M/G/1 system with
reservations (see \cite[p. 198]{BertGall92}) where the $n_s$
subregions correspond to users and the travel time between the
subregions corresponds to the reservation interval. Using the delay
expression for multiuser M/G/1 queue with reservations in \cite[p.
200]{BertGall92} we obtain,
\begin{equation}\label{eq:T_part}
T_{part}=\frac{\lambda s^2}{2(1-\rho)} +
\frac{n_s-\rho}{2v(1-\rho)}\sqrt{2}r^*+s,
\end{equation}
where $\rho=\lambda s$ is the system load. Combining this result
with (\ref{eq:T_LB_low_lam}) and noting that the above expression is
finite for all loads $\rho<1$, we have established the following
observation.
\begin{observation}
\emph{The time average delay in the system scales as
$\Theta(\frac{1}{1-\rho})$ with the load $\rho$ and the Partitioning
policy is stable for all $\rho <1$.}
\end{observation}
%

Despite the travel component of the service time, we can achieve
$\Theta(\frac{1}{1-\rho})$ delay as in classical queuing systems
(e.q., G/G/1 queue). This is the fundamental difference between this
system and the corresponding system where wireless transmission is
not used, as in the latter system the delay scaling with load is
$\Theta(\frac{1}{(1-\rho)^2})$ \cite{bert1}. This difference can be
explained intuitively as follows. Denote by
$N$ the average number of departures in a waiting time. 
It is easy to see from the P-K formula that in a classical M/G/1
queue, $N$ scales with the load as $\Theta(\frac{1}{1-\rho})$. We
argue that this scaling for $N$ is preserved in our system but not
in \cite{bert1}.
The $W_s$ expression as a function of $W_d$ in (\ref{eq:Ws_temp})
implies that for any given policy with its corresponding $W_d$, 
$N$ can be lower bounded by $\frac{\lambda W_d}{1-\rho}$. For the
system in \cite{bert1}, 
the minimum per-message distance the collector moves in the high
load regime scales as $\Omega(\frac{\sqrt{A}}{\sqrt{N}})$
\cite{bert1} (intuitively, the nearest neighbor distance among $N$
uniformly distributed points on a square region of area $A$ scales
as $\frac{\sqrt{A}}{\sqrt{N}}$). Therefore, for this system we have
$W_d \approx N\Omega(\frac{\sqrt{A}}{\sqrt{N}}) \approx
\Omega(\sqrt{NA})$ 
which gives $N \approx \Omega( \frac{\lambda^2 A}{(1-\rho)^2})$.
Namely, $W_d$ 
increases with the load and this results in an extra $1/(1-\rho)$
scaling in delay in addition to the $1/(1-\rho)$ factor of classical
M/G/1 queues. However, with the wireless reception capability, the
collector does not need to move for messages that are inside a disk
of radius $r^*$ around it. Since a finite (constant) number of such
disks cover the network region, $W_d$ can be upper bounded by a
constant independent of the system load (for the Partitioning policy
an easy upper bound on $W_d$ is the length of one cyclic tour around
the network). Therefore, in our system $N$ scales as $1/(1-\rho)$ as
in classical queues.

%

In \cite{Celik} we analyzed the case where the messages were
transmitted to the collector using a random access scheme (i.e., with probability $p$ in
each time slot) and
obtained $\Omega(\frac{1}{(1-\rho)^2})$ delay scaling as in the
system without wireless transmission. The reason for this is that in
order to have successful transmissions under the random access
interference of neighboring nodes, the reception distance should be
of the same order as the nearest neighbor distances \cite{Celik}, \cite{Tse}. 
%
%

\subsection{Numerical Results-Single Collector}\label{Sec:results-1}

\begin{figure}
\centering 
\includegraphics[width=0.4\textwidth]{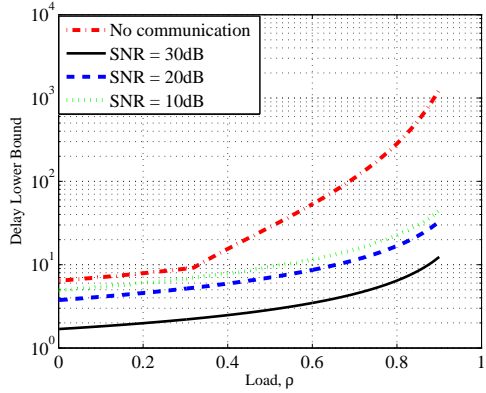}
\vspace{-3mm} \caption{Delay lower bound vs. network load using
different SNR values for $A=200$, $\beta=2$, $\alpha=4$, $v=1$ and
$s=1$.} \label{Fig:low_arr_delay_vs_lam_conts_r_SNRrangedB_noFad}
\vspace{-3mm}
\end{figure}

Here we present numerical results corresponding to the analysis in
the previous sections. We lower bound the delay expression in
(\ref{eq:T_LB_low_lam}) using $\mathbb{E} [ (||U||-r^*)^+] \ge
\mathbb{E}[||U||]-r^*$ (where $\mathbb[||U||]=0.383\sqrt{A}$ is
the expected distance of a uniform arrival to the center of square
region of area $A$ \cite{bert1}).
Fig.~\ref{Fig:low_arr_delay_vs_lam_conts_r_SNRrangedB_noFad} shows
the delay lower bound as a function of the network load for
different levels of channel quality\footnote{For the delay plot of
the no-communication system, the point that is not smooth arises
since the plot is the maximum of two delay lower bounds proposed in
\cite{bert1}.}. As the channel quality increases, the message delay
decreases as expected. For heavy loads, the delay in the system is
significantly less than the delay in the corresponding system
without wireless transmission in \cite{bert1}, demonstrating the
difference in the delay scaling between the two systems. For light
loads and more noisy communication channels, the delay performance
of the wireless network tends to the delay performance of
\cite{bert1}.

\begin{figure}
\centering 
\includegraphics[width=0.4\textwidth]{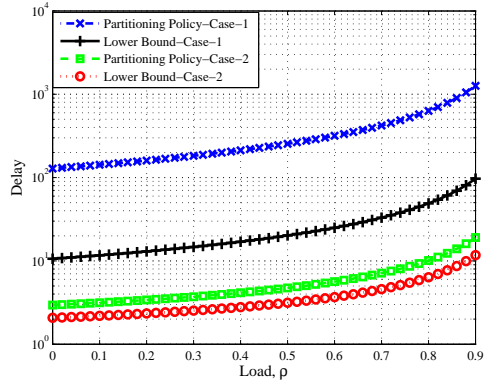}
\vspace{-3mm}\caption{Delay in the Partitioning policy vs the delay
lower bound for $SNR_c=17dB\, (r^*=2.2)$, $\beta=2$ and $\alpha=4$.
Case-1: Dominant travel time ($A=800$, $v=1$, $s=2$). Case-2:
Comparable travel and reception times ($A=60$, $v=10$, $s=2$).}
\label{Fig:partvsLB} \vspace{-4mm}
\end{figure}

Fig. \ref{Fig:partvsLB} compares the delay in the Partitioning
Policy to the delay lower bound for two different cases. When the
travel time dominates the reception time, the delay in the
Partitioning policy is about $10.6$ times the delay lower bound. For
a more balanced case, i.e., when the reception time is comparable to
the travel time, the delay ratio drops to $2.4$. 

\section{Multiple collectors}\label{Sec:Mult_veh}

In this section we extend our analysis to a wireless network with
multiple identical collectors. An arriving message is transmitted
when one of the $m$ collectors comes within the reception distance
of the message location and grants access for the message's
transmission. Therefore, at a given time there can be at most $m$
transmissions in the network. 
We consider policies that partition the network region into $m$
subregions. Each collector is assigned to one of the subregions and
is allowed to operate only in its own subregion. We call these
policies the class of \emph{network partitioning policies}. In such
a case, there is no interference from nodes within the subregion
where the transmission is taking place. The only source of
interference can be due to transmissions in other subregions.
However, we assume that signaling schemes used in different
subregions are orthogonal to each other so that there is no
destructive interference from other subregions. This can be achieved
for example by having different frequency bands for transmissions in
different subregions. Furthermore, for networks deployed on a large
area, even if orthogonal signaling is not utilized, the interference
between subregions can be negligible due to signal attenuation with
distance. Hence, we assume that simultaneous transmissions do not
interfere with each other. Note that this is consistent for the
purposes of a lower bound on delay since the interference from
neighboring nodes can only increase the message reception time.
Under these assumptions, utilizing the SNR criteria for successful
message reception as in Section \ref{Sec:LB_Delay}, the optimal
reception distance is given by $r^* = (\textrm{SNR}_c / \beta )^{1
/\alpha}$ and each reception takes time $s$.

\subsection{Stability}
Here we show that $\rho=\lambda s /m < 1$ is a necessary and sufficient condition for stability of the system.

\subsubsection{Necessary Condition for Stability}

A necessary condition for stability of the multi-collector system
is given by $\rho=\lambda s /m < 1$. 
We prove this by showing that the system stochastically dominates
the corresponding system with zero travel times (i.e., an M/D/m
queue, a queue with Poisson arrivals, constant service time and $m$
servers) similar to Section \ref{Sec:Stab}.
\begin{theorem}\label{thm:mt_LB_delay1}
\emph{A necessary condition for the stability of any policy is $\rho=\lambda s /m < 1$.
Furthermore, 
the optimal steady state time average delay $T_m^*$ is lower
bounded by
\begin{equation}\label{eq:mt_W_temp1}
T_m^* \ge \frac{\lambda s^2}{2m^2(1-\rho)}  -
\frac{m-1}{m}\frac{s^2}{2s} +s,
\end{equation}
where $\rho = \lambda s/m$ is the system load.}
\end{theorem}
The proof is similar to the proof of Theorem \ref{thm:inf_vel_sin_veh} and is given in Appendix-A. It makes use of the fact that the steady state time average delay in the system is at least
as big as the delay in the equivalent system in which travel times
are considered to be zero (i.e., $v=\infty$).

\subsubsection{Sufficient Condition for Stability}

We argue that $\rho<1$ is also sufficient for stability. This can be seen by dividing the network region into $m$ equal
subregions, and performing a single-collector TSPN policy in each
subregion. Since the arrival process is Poisson, each subregion receives an independent Poisson arrival process of intensity $\lambda/m$. Furthermore, each collector performs a  TSPN policy independently of the other collectors. Therefore, using the stability result of the single-collector TSPN policy, the systems in each subregion are stable if $\rho<1$. We state this fact in the following theorem:
\begin{theorem}\label{thm:mt_LB_delay1}
\emph{The system is stable under the multi-collector TSPN policy for all loads $\rho=\lambda s / m < 1$.}
\end{theorem}
Note that a similar delay analysis to the single-collector TSPN case shows that the multi-collector TSPN policy has $O(\frac{1}{1-\rho})$ delay scaling with the load $\rho$.

\subsection{Delay Lower Bound}
In addition to the delay lower bound given in
(\ref{eq:mt_W_temp1}), we provide another lower
bound for the optimal delay $T_m^*$ and take their simple average. The following lemma 
states the second lower bound on delay. It is based on the
convexity argument that when travel component of the waiting time
is lower bounded by a constant, 
the equal area partitioning of the
network region minimizes the resulting delay expression out of all
area partitionings. 
\begin{theorem}\label{thm:mt_LB_delay2}
\emph{For the class of network partitioning policies, the optimal
steady state time average delay $T_m^*$ is lower bounded by}
\begin{equation}\label{eq:mt_W_temp2}
T_m^* \ge \frac{1}{1-\rho}
\frac{\textrm{max}\big(0,\frac{2}{3}\sqrt{\frac{A}{m \pi}}-r^*\big)
}{v} +s,
\end{equation}
\emph{where $\rho = \lambda s/m$ is the system load.}
\end{theorem}
%
%
%
\begin{proof}
Here we use an approach similar to the proof of Theorem
\ref{thm:LB_delay}. We divide the average delay
$T$ into three components: 
\begin{equation}\label{eq:m_T_LB_temp}
T=W_{d} + W_{s} + s.
\end{equation}
The below lemma provides a bound for $W_d$, the average message
waiting time due to the collectors's travel, using a result in
\cite{Haim} for the \emph{m}-median problem.
\begin{lemma}\label{lem:app_pf2}
\begin{equation}\label{eq:m_T_LB_temp3}
W_{d} \ge \frac{\textrm{max}\big(0,\frac{2}{3}\sqrt{\frac{A}{m
\pi}}-r^*\big) }{v}.
\end{equation}
\end{lemma}
\begin{proof}
Let $\Omega$ be \emph{any} set of points in $\Re$ with $|\Omega|=m$.
Let $U$ be a uniformly distributed location in $\Re$ independent of
$\Omega$ and define $Z^* \triangleq \textrm{min}_{\nu \in \Omega}
\parallel U-\nu \parallel$.
%
Let the random variable $Y$ be the distance from the center of a
disk of area $A/m$ to a uniformly distributed point within the disk.
Then it is shown in \cite{Haim} that 
\begin{equation}\label{eq:haim_bound}
\mathbb{E}[f(Z^*)] \ge \mathbb{E}[f(Y)]
\end{equation}
for any nondecreasing function $f(.)$.
Using this result we obtain $\mathbb{E}[\max(0,Z^*-r^*)] \ge
\mathbb{E}[\max(0,Y-r^*)]$.
Note that $W_d$ can be lower bounded by the expected distance of a
uniform arrival to the closest collector at the time of arrival less
$r^*$. Because the travel distance is nonnegative, we have
\begin{equation*}
W_{d} \ge \mathbb{E}[\max(0,Y-r^*)]/v \ge
\textrm{max}(0,\mathbb{E}[Y]-r^*)/v,
\end{equation*}
where the second bound is due to Jensen's inequality. Substituting
$\mathbb{E}[Y]=\frac{2}{3}\sqrt{\frac{A}{m \pi}}$ into the above
expression completes the proof.
\end{proof}
Intuitively the best a priori placement of $m$ points in $\Re$ in
order to minimize the distance of a uniformly distributed point in
the region to the closest of these points is to cover the region
with $m$ disjoint disks of area $A/m$ and place the points at the
centers of the disks. Such a partitioning of the region is not
possible, however, using this idea we can lower bound the expected
distance as in (\ref{eq:haim_bound}).
%

We now derive a lower bound on $W_s$. Let $R^1,R^2,...,R^m$ be the
network partitioning with areas $A^1,A^2,...,A^m$ respectively
($\sum_{j=1}^{m}A^j=A$).
Consider the message receptions in steady state 
that are received by collector $j$ eventually. Let $\lambda^j$ be
the fraction of the arrival rate served by collector $j$. Due to the
uniform distribution of the message locations we have
\begin{equation*}
\frac{\lambda^j}{\lambda}=\frac{A^j}{A}.
\end{equation*}
Let $N^j$ be the average number of message receptions for which the
messages that are served by collector $j$ waits in steady state.
Similarly let $W^j_s$ and $W^j_d$ be the average waiting times for
messages served by collector $j$ due to the time spent on message
receptions and collector $j$'s travel respectively. Using
(\ref{eq:Ws_temp1}) and lower bounding the residual time by zero we
have
\begin{equation*}
W_s^j \ge s N^j.
\end{equation*}
Using Little's law ($N^j=\lambda^j (W_s^j+W_d^j)$) similar to the
derivation of (\ref{eq:Ws_temp}) we have
\begin{equation}\label{eq:Ws^j_temp}
W_s^j \ge \frac{\lambda^j s}{1-\lambda^j s}W_d^j.
\end{equation}
The fraction of messages served by collector $j$ is $A^j/A$.
Therefore, we can write $W_s$ as
\begin{eqnarray}
W_s  \!\!\!\!\!\!\!\!\!\!&& = \displaystyle \sum_{j=1}^m \frac{A^j}{A}W_s^j \nonumber\\
&& \ge  \displaystyle \sum_{j=1}^m \frac{A^j}{A}\frac{\lambda^j
s}{1-\lambda^j s}W_d^j.\label{eq:Ws^j_temp4}
\end{eqnarray}
For a given region $R^j$ with area $A^j$, $W_d^j$ is lower bounded
by (similar to the derivation of (\ref{eq:Wd_temp})) the distance of
a uniform arrival to the median of the region less $r^*$.
\begin{eqnarray}
W_d^j \!\!\!\!\!\!\!\!\!\!&& \ge \frac{\mathbb{E} [ \max(0,||U-\nu||-r^*)] }{v}\nonumber\\
&& \ge \frac{\max(0,\mathbb{E}
[||U-\nu||]-r^*)}{v},\label{eq:Ws^j_temp5}
\end{eqnarray}
where $\nu$ is the median of $R^j$ and $||U-\nu||$ is the distance
of $U$, a uniformly distributed location inside $R^j$, to $\nu$. The
inequality in (\ref{eq:Ws^j_temp5}) is due to Jensen's inequality
for convex functions. A disk shaped region yields the minimum
expected distance of a uniform arrival to the median of the region.
Using this we further lower bound $W_d$ by noting that for a disk
shaped region of area $A_j$, $\mathbb{E} [||U-\nu||]$ is just the
expected distance of a uniform arrival to the center of the disk
given by $\frac{2}{3}\sqrt{\frac{A^j}{\pi}}$. Hence
\begin{equation}\label{eq:Ws^j_temp3}
W_d ^j\ge
\frac{\textrm{max}(0,\frac{2}{3}\sqrt{\frac{A^j}{\pi}}-r^*)}{v} =
\frac{\textrm{max}(0,c_1 \sqrt{A^j}-r^*)}{v},
\end{equation}
where $c_1=\frac{2}{3\sqrt{\pi}}=0.376$. Letting
$f(A^j)=\frac{\lambda \frac{A^j}{A} s}{1-\lambda \frac{A^j}{A} s}$,
which is a convex and increasing function of $A^j$, we rewrite
(\ref{eq:Ws^j_temp4}) as
\begin{equation}\label{eq:Ws^j_temp6}
W_s \ge \displaystyle \sum_{j=1}^m \frac{f(A^j)}{vA} A^j
\textrm{max}(0,c_1 \sqrt{A^j}-r^*).
\end{equation}
Next we will show that the function $f(A^j)A^j \textrm{max}(0,c_1
\sqrt{A^j}-r^*)$ is a convex function of $A^j$ via the two lemmas
below.
\begin{lemma}\label{lem:app_pf3}
\emph{Let $f(.)$ and $g(.)$ be two convex and increasing functions
(possibly nonlinear) defined on $[0,A]$. The function $h(.)=f.g(.)$
is also convex and increasing on its domain $[0,A]\textrm{x}[0,A]$.}
\end{lemma}
\begin{proof}
See Appendix-B.
\end{proof}

\begin{lemma}\label{lem:app_pf4}
\emph{$h(x)=x \max(0,c_1\sqrt{x}-c_2)$ is a convex and
increasing function of $x$.}
\end{lemma}
\begin{proof}
See Appendix-C.
\end{proof}
Letting $g(A^j) \doteq f(A^j)A^j \textrm{max}(0,c_1
\sqrt{A^j}-r^*)$, we have from the lemmas \ref{lem:app_pf3} and
\ref{lem:app_pf4} that the function $g(A^j)$ is convex. Now
rewriting (\ref{eq:Ws^j_temp6}) we have
\begin{equation*}
W_s \ge \displaystyle (\frac{m}{vA})\frac{1}{m}\sum_{j=1}^m g(A^j).
\end{equation*}
Using the convexity of the function $g(A^j)$ we have
\begin{eqnarray}
W_s \!\!\!\!\!\!\!\!\!\!&& \ge \displaystyle
(\frac{m}{vA})g\Big(\frac{\sum_{j=1}^m
A^j}{m} \Big) \nonumber\\
&&=\frac{m}{vA}g(\frac{A}{m})  \nonumber\\
&&= \frac{\frac{\lambda s}{m}}{1-\frac{\lambda s}{m}}
\frac{\textrm{max}(0,c_1
\sqrt{A/m}-r^*)}{v} \nonumber\\
&& = \frac{\rho}{1-\rho} \frac{\textrm{max}\big(c_1
\sqrt{\frac{A}{m}}-r^*\big)}{v}. \label{eq:Ws^j_temp7}
\end{eqnarray}
The above analysis essentially implies that the $W_s$ expression in
(\ref{eq:Ws^j_temp6}) is minimized by the \emph{equitable
partitioning} of the network region. Finally combining
(\ref{eq:m_T_LB_temp}), (\ref{eq:m_T_LB_temp3}) and
(\ref{eq:Ws^j_temp7}) we obtain (\ref{eq:mt_W_temp2}).
%
\end{proof}

Finally, taking the simple average of (\ref{eq:mt_W_temp2}) and
(\ref{eq:mt_W_temp1}) we arrive at the following theorem.
\begin{theorem}\label{thm:mt_LB_delay}
\emph{For the class of network partitioning policies, the optimal
steady state time average delay $T_m^*$ is lower bounded by}
\begin{equation}
T_m^* \ge \! \frac{\lambda s^2}{4m^2(1\!-\!\rho)} +\!
\frac{\textrm{max}\big(0,\frac{2}{3}\sqrt{\frac{A}{m
\pi}}\!-\!r^*\big) }{2v(1\!-\!\rho)}
-\!\frac{m\!-\!1}{m}\frac{s^2}{4s}+ s,\label{eq:Tm*_highload}
\end{equation}
\emph{where $\rho = \lambda s/m$ is the system load.}
\end{theorem}
Theorem \ref{thm:mt_LB_delay} is valid for the class of network
partitioning policies. For the system without wireless transmission,
it has been shown that partitioning the region into $m$ equal size
disjoint subregions (one for each collector) preserves optimality in
the high load limit \cite{bert2}, \cite{Xu}. We conjecture that this
optimality is also preserved in our system.

\subsection{Multiple Collector Policies}
Note that it is shown in \cite{bert2} that a generalization of the
FCFS policy, namely, creating $m$ Voronoi regions with centers of
the regions given by the $m$-median locations\footnote{The set of
$m$-median locations for a region is the set of the best $m$ a
priori locations in the region that minimizes the expected distance
to a uniform arrival.} of the network region and having each
receiver perform the single-receiver FCFS service in each region has
optimal delay at light loads. This also holds in our system via a
similar argument to Section \ref{Sec:Policies}. Moreover, this
policy is not stable as $\rho\rightarrow 1$ due to the same reason
as in Section \ref{Sec:Policies}.

\subsubsection{Generalized Partitioning Policy}
\begin{figure}
\centering \psfrag{1\r}[l][][.5]{$\!r^*$}
\psfrag{2\r}[l][][.5]{$\!\!\!\!\!\sqrt{2}r^*$}
\includegraphics[width=0.25\textwidth]{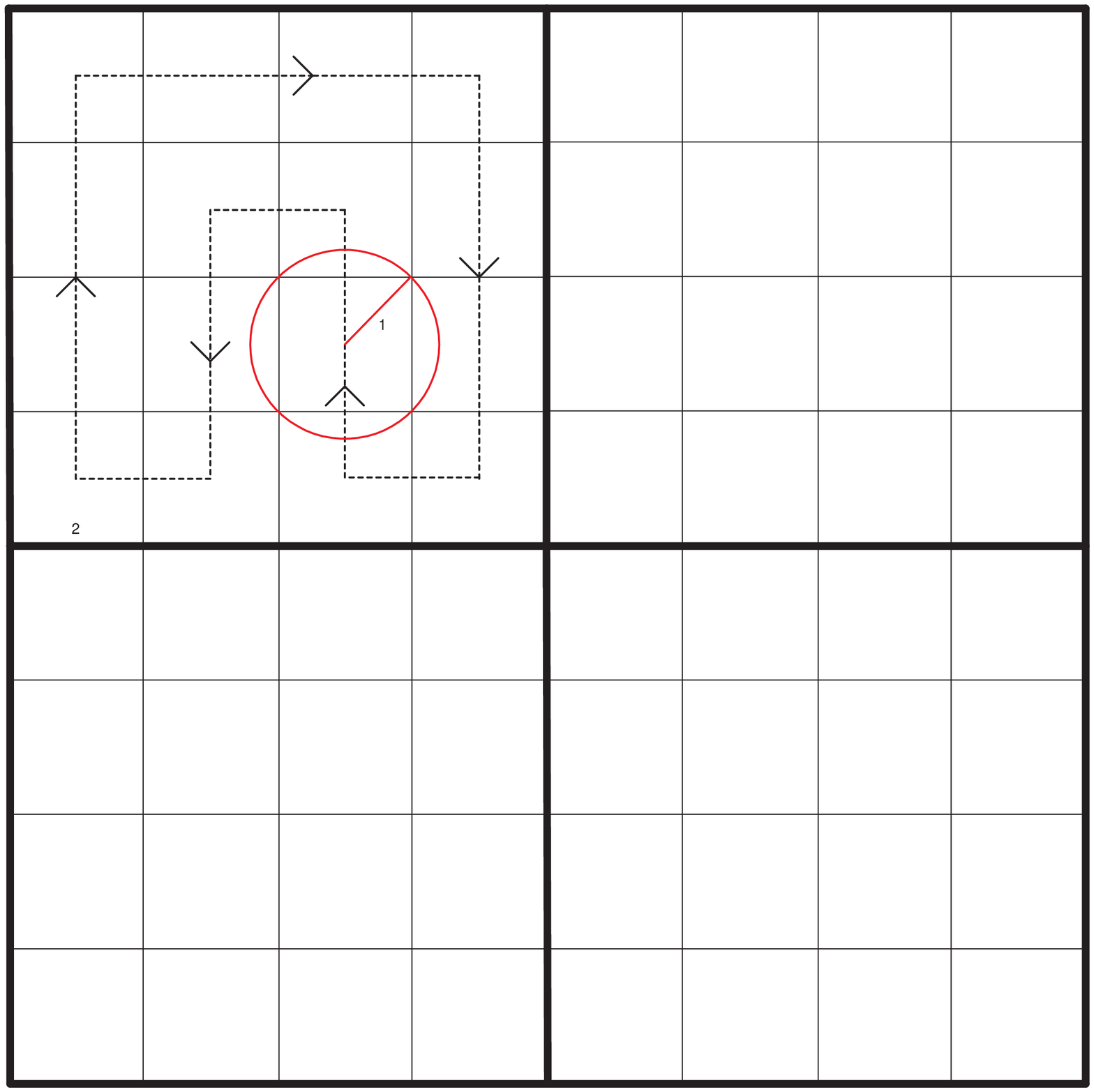}
\vspace{-1mm}\caption{The partitioning of the network region into
square subregions of side $\sqrt{2}r^*$ for the case of multiple
collectors in the network. The circle with radius $r^*$ represents
the communication range and the dashed lines represent one of the
collector's path.} \label{Fig:part_pol_m}
\end{figure}

Next we propose a policy based on dividing the network region into
$m$ equal size subregions. Each collector is assigned to one of the
subregions and is responsible for receiving messages that arrive
into its own subregion using the single collector partitioning
policy analyzed in Section \ref{Sec:Part_pol}. Namely, first the
network region is divided into subregions of area $A/m$ and then
each subregion is divided into
$\sqrt{2}r^* \,\textrm{x}\, \sqrt{2}r^*$ squares\footnote{Note that
the number of collectors must be a square number in order to divide
the network region into subregions of exactly equal shape and
size.}. The number of $\sqrt{2}r^* \,\textrm{x}\, \sqrt{2}r^*$
squares in each subregion is given by $n_s=\frac{A/m}{2(r^*)^2}$.
Fig. \ref{Fig:part_pol_m} represents such a partitioning for the
case of four collectors in the network with $n_s=16$ squares in each
subregion.
 %
%
%
%
%
Since each subregion behaves identically, the average delay of this
policy is the average delay of the single collector Partitioning
policy applied to a subregion with arrival rate $\lambda/m$, area
$A/m$, and $n_s=\frac{A/m}{2(r^*)^2}$:
\begin{equation}\label{eq:mT_part}
T_{part}=\frac{\lambda s^2}{2m(1-\rho)} +
\frac{\frac{A}{2m(r^*)^2}-\rho}{2v(1-\rho)}\sqrt{2}r^*+s,
\end{equation}
where $\rho=\lambda s/m$ is the system load. This result, when
combined with (\ref{eq:mt_W_temp1}), establishes that for the case
of multiple collectors in the system the delay scaling with the load
is $\Theta(\frac{1}{1-\rho})$. This is again a fundamental
improvement compared to the $\Theta(\frac{1}{(1-\rho)^2})$ delay
scaling in the system without wireless transmission and with
multiple collectors in \cite{bert2}.

\subsubsection{Numerical Results}
\begin{figure}
\centering 
\includegraphics[width=0.4\textwidth]{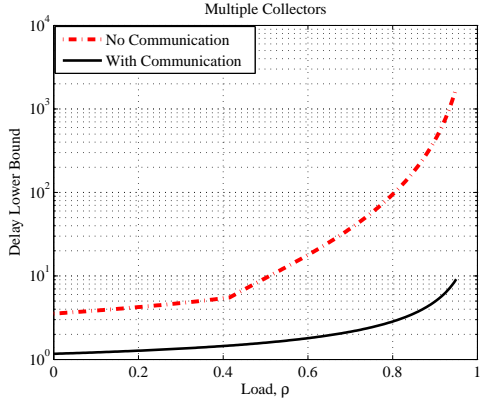}
\vspace{-3mm}\caption{Delay lower bound vs. network load for m=2
collectors, $SNR_c=30dB\,(r^*=4.7)$, $A=400$, $\beta=2$, $\alpha=4$,
$v=1$ and $s=1$.}\label{Fig:singlevsmultivsnocomm} \vspace{-4mm}
\end{figure}
\begin{figure}
\centering 
\includegraphics[width=0.4\textwidth]{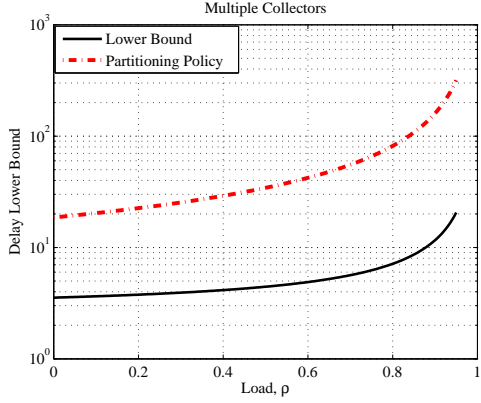}
\vspace{-3mm}\caption{Delay of the Partitioning policy vs the delay
lower bound for $m=4$ collectors, $SNR_c=20dB\, (r^*=2.6)$, $A=500$,
$\beta=2$, $\alpha=4$, $v=1$ and $s=2$.}\label{Fig:partmvsLBm}
\vspace{-4mm}
\end{figure}

We compare the delay lower bound in (\ref{eq:Tm*_highload}) to the
delay lower bound in the corresponding system without wireless
transmission in \cite {bert2} for the case of two collectors in
Fig.~\ref{Fig:singlevsmultivsnocomm}. The delay in the two-collector
system is significantly below the delay in the system without
wireless transmission and this difference is more pronounced for
high loads. 
Fig.~\ref{Fig:partmvsLBm} displays the delay lower bound in
(\ref{eq:Tm*_highload}) and the delay of the Partitioning policy in
(\ref{eq:mT_part}) as functions of the network load $\rho$. The
delay of the Partitioning policy is about $7$ times the delay lower
bound. 

\section{CONCLUSION}\label{Sec:Conc}
In this paper we considered the use of dynamic vehicle routing in
order to improve the delay performance of 
wireless networks where messages arriving randomly in time and space
are gathered by a mobile collector.
We characterized the stability region of this system to be all
system loads $\rho <1$ and derived fundamental lower bounds on time
average expected delay. We derived upper bounds on delay by
analyzing policies and 
extended our results to the case of multiple collectors in the
system. Our results show that combining controlled mobility and
wireless transmission results in $\Theta(\frac{1}{1-\rho})$ delay scaling with load $\rho$. 
This is the fundamental difference between our system and the system
without wireless transmission (DTRP) analyzed in \cite{bert1} and
\cite{bert2} where the delay scaling with the load is
$\Theta(\frac{1}{(1-\rho)^2})$.

This work is a first attempt towards utilizing a combination of
controlled mobility and wireless transmission for data collection in
stochastic and dynamic wireless networks. Therefore, there are many
related open problems. In this paper we have utilized a simple
wireless communication model based on a communication range. In the
future we intend to study more advanced wireless communication
models such as modeling the transmission rate as a function of the
transmission distance. 
Finally, extending our results for a general message
location distribution in the network is a subject of future
research.


%


\section*{Appendix A-Proof of Theorem \ref{thm:mt_LB_delay1}}

The proof is similar to the proof of Theorem \ref{thm:inf_vel_sin_veh}. First consider the following lemma.
\begin{lemma}\label{lem:app_pf}
\emph{The steady state time average delay in the system is at least
as big as the delay in the equivalent system in which travel times
are considered to be zero (i.e., $v=\infty$).}
\end{lemma}
%
%
\begin{proof}
The proof is similar to the proof of Theorem \ref{lem:app_pf}.
Consider the summation of per-message reception and travel times,
$s$ and $d_i$, as the total service requirement of a message in each
system. Since $d_i$ is zero for all $i$ in the infinite velocity
system and since the reception times are constant equal to $s$ for
both systems, the total service requirement of each message in our
system is deterministically greater than that of the same message in
the infinite velocity system. Let $D_1,D_2,...,D_n$ be the departure
instants of the first, second and similarly the $n^{th}$ message in
the system. Similarly let $D^{'}_1,D^{'}_2,...,D^{'}_n$ be the
departure instants of the first, second and similarly the $n^{th}$
message in the infinite velocity system. Similarly let
$A_1,A_2,...,A_n$ be the arrival times of the first second and
similarly the $n^{th}$ message in both systems. We will use
\emph{complete induction} to prove that $D_i \ge D_i^{'}$ for all
$i$. Since the service requirement of each message is less in the
infinite velocity system, we have $D_1 \ge D_1^{'}$. Assume we have
$D_i \ge D_i^{'}$ for all $i\le n$. We need to show that $D_{n+1}
\ge D_{n+1}^{'}$ in order to complete the complete induction. We
have
\begin{equation}\label{eq:app_proof_temp1}
A_{n+1} \le D_{n+1}-s,
\end{equation}
hence the $n+1^{th}$ message is available at time $D_{n+1}-s$. We
also have
\begin{equation*}
D_{n+1-m}^{'} \le D_{n+1-m} \le D_{n+1}-s.
\end{equation*}
The first inequality is due to the complete induction hypothesis and
the second inequality is due the fact that the $m^{th}$ last
departure before the $n+1^{th}$ departure has to occur before the
time $D_{n+1}-s$. Hence there is at least one collector available in
the infinite velocity system before the time $D_{n+1}-s$. Combining
this with (\ref{eq:app_proof_temp1}) proves the complete induction.
Now let $D(t)$ and $D^{'}(t)$ be the total number of departures by
time $t$ in our system and the infinite velocity system
respectively. Similarly let $N(t)$ and $N^{'}(t)$ be the total
number of messages in the two systems at time $t$. Finally let
$A(t)$ be the total number of arrivals by time $t$ in both systems.
We have $N(t)=A(t)-D(t)$ and $N^{'}(t)=A(t)-D^{'}(t)$. From the
above induction we have $D(t)\le D^{'}(t)$ and therefore $N(t) \ge
N^{'}(t)$. Since this is true at all times, we have that the time
average number of customers in the system is greater than that in
the infinite velocity system. Finally using Little's law proves the
lemma.
\end{proof}

When the travel time is considered to be zero, the system becomes an
M/D/m queue (a queue with Poisson arrivals, constant service time
and $m$ servers). Therefore we can bound $T_m^*$ using bounds for
general G/G/m systems. In particular, the waiting time $W_{G/G/m}$
in a G/G/m queue with service time $s$ is bounded below by \cite[p.
48]{Klein76}
\begin{equation}\label{eq:mt_W_temp}
W_{G/G/m} \ge \hat{W} - \frac{m-1}{m}\frac{s^2}{2s},
\end{equation}
where $\hat{W}$ is the waiting time in a single server system with
the same arrivals as in the G/G/m queue and service time $s/m$.
Since in our case the infinite velocity system behaves as an M/D/m
system, $\hat{W}$ has an exact expression given by the P-K formula:
$\hat{W}=\lambda s^2/(2m^2(1-\rho))$ where $\rho = \lambda s/m$.
Substituting this in (\ref{eq:mt_W_temp}) and using Lemma
\ref{lem:app_pf} we have (\ref{eq:mt_W_temp1}).

\section*{Appendix B-Proof of Lemma \ref{lem:app_pf3}}

Clearly $h(.)=f.g(.)$ is increasing. Let $x$ and $y$ be two points in the
domain of $h$ and let $\alpha \in (0,1)$ be a real number.
\begin{eqnarray}
h(\alpha x + (1-\alpha) y) \!\!\!\!\!\!\!\!\!\!&& = f(\alpha x + (1-\alpha) y)g(\alpha x + (1-\alpha) y) \nonumber\\
&& \le (\alpha f(x) + (1-\alpha)f(y)).\nonumber\\
&& \;\;\;\, .(\alpha g(x)+(1-\alpha)g(y))\nonumber\\
&& = \alpha^2 f(x)g(x) + (1-\alpha)^2f(y)g(y) \nonumber\\
&& \;\;\;\, +
\alpha(1\!-\!\alpha)f(x)g(y)\!+\!\alpha(1\!-\!\alpha)f(y)g(x)),\nonumber
\end{eqnarray}
where the inequality is due to the convexity of $f$ and $g$. We add
and subtract $\alpha f(x)g(x)$ and after some algebra obtain
\begin{eqnarray}
h(\alpha x + (1-\alpha) y) \!\!\!\!\!\!\!\!\!\!&& \le \alpha h(x) +(1-\alpha)h(y) \nonumber\\
&&\;\;\;\, + \alpha(1-\alpha)(f(x)-f(y))(g(y)-g(x))  \nonumber\\
&& \le \alpha h(x) + (1-\alpha)h(y), \nonumber
\end{eqnarray}
where the last inequality is due to the fact that $f$ and $g$ are
increasing functions.

\section*{Appendix C-Proof of Lemma \ref{lem:app_pf4}}

It is clear that $h(x)$ is an increasing function of $x$. Let $x$
and $y$ be two points in the domain of $h$ and let $\alpha \in
(0,1)$ be a real number.
\begin{eqnarray}
h(\alpha x \!\!\!\!\!\!\!\!\!\!\!&&+  (1-\alpha) y) = \nonumber\\
&& =(\alpha x + (1-\alpha) y)\textrm{max}(0,c_1\sqrt{\alpha x + (1-\alpha) y}-c_2) \nonumber\\
&& = \textrm{max}(0,c_1(\alpha x + (1-\alpha) y)^{\frac{3}{2}}-c_2(\alpha x + (1-\alpha) y)) \nonumber\\
&& \le \textrm{max}\big(0,c_1(\alpha x^{\frac{3}{2}} + (1-\alpha) y^{\frac{3}{2}})-c_2(\alpha x + (1-\alpha) y)\big) \nonumber\\
&& = \textrm{max}\big(0,\alpha x (c_1\sqrt{x}-c_2) + (1-\alpha) y (c_1\sqrt{y}-c_2)  \big)  \nonumber\\
&& \le \textrm{max}\big(0,\alpha x (c_1\sqrt{x}\!-\!c_2)\big) \!+\! \textrm{max}\big(0,(1\!-\!\alpha) y (c_1\sqrt{y}\!-\!c_2)  \big)  \nonumber\\
&& = \alpha h(x) + (1-\alpha) h(y),  \nonumber
\end{eqnarray}
where the first inequality is due to the convexity of the function
$x^{\frac{3}{2}}$. This shows that $h(x)$ is a convex and increasing function.


%
%

%


\bibliographystyle{IEEEtranS}

\begin{thebibliography}{99}

\bibitem{Akyildiz_05} I. F. Akyildiz, D. Pompili, and T. Melodia, ``Underwater Acoustic
Sensor Networks: Research Challenges,'' \emph{Ad Hoc Networks
(Elsevier)}, vol.~3, no.~3, pp.~257-279, Mar. 2005.

\bibitem{AltKonsLiu92} E. Altman, P. Konstantopoulos, and Z. Liu, ``Stability, monotonicity and invariant quantities in
general polling systems,'' \emph{Queuing Sys.}, vol.~11, pp.~35-57,
1992.
%

\bibitem{Asmussen} S. Asmussen, ``Applied Probability and Queues,'' 
Wiley, 1987.

\bibitem{arkinGCSP94} E. M. Arkin and R. Hassin, ``Approximation algorithms for the geometric covering salesman problem,''
\emph{Discrete Applied Mathematics}, vol.~55, pp.~197-218, 1994.




\bibitem{BertGall92} D. Bertsekas and R. Gallager, ``Data Networks,'' Prentice Hall, 1992.


\bibitem{bert1} D. J. Bertsimas and G. van Ryzin, ``A stochastic and dynamic vehicle routing problem in the Euclidean plane,''
\emph{Opns. Res.}, vol.~39, pp.~601-615, 1990.

\bibitem{bert2} D. J. Bertsimas and G. van Ryzin, ``Stochastic and dynamic vehicle routing in the Euclidean
plane with multiple capacitated vehicles,'' \emph{Opns. Res.},
vol.~41, pp.~60-76, 1993.

\bibitem{bert3} D. J. Bertsimas and G. van Ryzin, ``Stochastic and dynamic vehicle
routing with general demand and interarrival time distributions,''
\emph{Adv. App. Prob.}, vol.~20, pp.~947-978, 1993.

\bibitem{Burns-Levine-Inf05} B. Burns, O. Brock, and B. N. Levine, ``MV routing and capacity building in disruption tolerant networks,'' In \emph{Proc. IEEE INFOCOM'05}, Mar. 2005.

\bibitem{Celik} G. D. Celik and E. Modiano, ``Random access wireless networks with controlled mobility'', In \emph{Proc. IFIP MEDHOCNET'09}, Jun. 2009.




%



\bibitem{FB2004} E. Frazzoli and F. Bullo, ``Decentralized algorithms for vehicle routing in a stochastic time-varying environment,''
In \emph{Proc. IEEE CDC'04}, Dec. 2004.

\bibitem{Geor_Neely_Tass06} L. Georgiadis, M. Neely, and L. Tassiulas, ``Resource Allocation and Cross-Layer Control in Wireless Networks,'' Now Publishers, 2006.



\bibitem{GMPS04}{ A. E. Gammal, J. Mammen, B. Prabhakar, and D. Shah, ``Throughput-delay trade-off in wireless networks,'' In \emph{Proc. IEEE INFOCOM'04}, Mar.
2004.}



\bibitem{Tse} M. Grossglauser and D. Tse, ``Mobility increases the capacity of ad hoc wireless networks ,'' \emph{IEEE/ACM Trans.
Netw.}, vol. 11, no. 1, pp. 125–-137, Feb. 2003.



\bibitem{GuptaKumar} P. Gupta and P. R. Kumar, ``The capacity of wireless networks,'' \emph{IEEE Trans. Inf. Theory}, vol. 46, no. 2, pp.~388–-404, Mar. 2000.


\bibitem{Haim} M. Haimovich and T. L. Magnanti, ``Extremum properties of of hexagonal partitioning and the uniform distribution in euclidian location,'' \emph{SIAM J. Disc. Math.} 1, 50-64, 1988.

\bibitem{Ver-Altm-09} V. Kavitha and E. Altman, ``Queueing in Space:
design of Message Ferry Routes in sensor networks,'' In \emph{Proc. ITC'09}, Sep. 2009.

%





\bibitem{Klein76} L. Kleinrock, ``Queueing Systems: Volume 2: Computer Applications,'' John Wiley $\&$ Sons, 1976.

\bibitem{Lawler85} E. L. Lawler, J.Lenstra, A. Kan, D. Shmoys, `` The Traveling salesman problem :a guided
tour of combinatorial optimization,'' Wiley, 1985.



%


\bibitem{Luo-Hub05} J. Luo and J. P. Hubaux, ``Joint mobility and routing for lifetime elongation in wireless sensor networks,'' In \emph{Proc. IEEE INFOCOM'05}, Mar. 2005.

\bibitem{MitTSPN07} J. S. B. Mitchell, ``A PTAS for TSP with neighborhoods among fat regions in the plane,''
In \emph{Proc. ACM-SIAM SODA'07}, Jan. 2007.

\bibitem{ModShahZuss06} E. Modiano, D. Shah and G. Zussman, ``Maximizing throughput in
wireless networks via Gossip,'' In \emph{Proc. ACM
SIGMETRICS/Performance'06}, June 2006.



\bibitem{NM05} M. J. Neely and E. Modiano, `` Capacity and delay tradeoffs for ad hoc mobile networks,'' \emph{IEEE Trans. Inf. Theory}, vol. 51, no. 6, pp. 1917–-1937, Jun. 2005.

%
\bibitem{PavFraz07} M. Pavone, N. Bisnik, E. Frazzoli and V. Isler, ``Decentralized vehicle routing in a stochastic and dynamic environment
with customer impatience,'' In \emph{Proc. RoboComm'07}, Oct. 2007.

\bibitem{Jerome08} J. Le Ny, M. Dahleh, E. Feron, and E. Frazzoli, ``Continuous Path
Planning for a Data Harvesting Mobile Server,'' In \emph{Proc. IEEE
CDC'08}, Dec. 2008.


\bibitem{emilio} V. Sharma, E. Frazzoli, and P. G. Voulgaris, ``Delay in mobility-assisted constant-throughput wireless networks,'' In \emph{Proc. IEEE CDC'05}, Dec. 2005.

\bibitem{SH08} Y. Shi and Y. T. Hou, ``Theoretical results on base station movement problem for sensor network,'' In \emph{Proc. IEEE INFOCOM'08}, Apr.
2008.

\bibitem{Soma-kansal-06} A. Somasundara, A. Kansal, D. Jea, D. Estrin, and M. Srivastava, ``Controllably mobile infrastructure for low energy
embedded networks,'' \emph{IEEE Trans. Mobile Comput.}, vol. 5, no.
8, pp. 958-973, Aug. 2006.


\bibitem{Holly06} H. Waisanen, D. Shah, M. A. Dahleh, `` Minimal Delay in Controlled Mobile Relay Networks,'' In \emph{Proc. IEEE CDC'06}, Dec. 2006.

\bibitem{Holly07} H. Waisanen, D. Shah, M. A. Dahleh, `` Lower Bounds for Multi-Stage Vehicle Routing,'' In \emph{Proc. IEEE CDC'07}, Dec. 2007.


\bibitem{Xu} H. Xu, ``Optimal policies for stochastic and dynamic vehicle routing problems,'' Ph.D. Thesis, MIT, 1994.


\bibitem{ZhaMaYa08} M. Zhao, M. Ma, and Y. Yang, ``Mobile data gathering with space-division multiple access in wireless sensor
networks,'' In \emph{Proc. IEEE INFOCOM'08}, Apr. 2008.


\bibitem{ZAZ2004} W. Zhao, M. Ammar, and E. Zegura, ``A message ferrying approach for data delivery in sparse mobile ad hoc networks,'' In \emph{Proc. ACM MobiHoc'04}, May 2004.




\end{thebibliography}

\end{document}